\documentclass[11pt]{article}

\usepackage[utf8]{inputenc}
\usepackage[T1]{fontenc}
\usepackage{lmodern}
\usepackage{amsmath,amssymb,amsthm,mathtools}
\usepackage[a4paper,margin=1in]{geometry}
\usepackage{hyperref}
\usepackage{enumitem}
\usepackage{microtype}
\usepackage{mathrsfs}

\hypersetup{colorlinks=true,linkcolor=blue,citecolor=blue,urlcolor=blue}

\newtheorem{theorem}{Theorem}[section]
\newtheorem{lemma}[theorem]{Lemma}

\newtheorem{corollary}[theorem]{Corollary}
\theoremstyle{definition}
\newtheorem{definition}[theorem]{Definition}
\theoremstyle{remark}
\newtheorem{remark}[theorem]{Remark}
\newcommand{\rk}{\operatorname{rank}}
\newcommand{\Bin}{\operatorname{Bin}}

\newcommand{\rank}{\operatorname{rank}}
\newcommand{\floor}[1]{\left\lfloor #1 \right\rfloor}

\title{Majority Edge-Colourings of Hypergraphs}
\author{
Jiangdong Ai\thanks{School of Mathematical Sciences and LPMC, Nankai University. {\tt jd@nankai.edu.cn}. Funded by the National Natural Science Foundation of China (No.12522117, No.12401456), the Natural Science Foundation of Tianjin (No.24JCQNJC01960) and Fundamental and Interdisciplinary Disciplines Breakthrough Plan of the Ministry of Education of China (JYB2025XDXM207).}
\and
Feiyu Nan\thanks{School of Mathematical Sciences and LPMC, Nankai University, Tianjin 300071, P.R. China. Email: \texttt{2211185@mail.nankai.edu.cn}.}
}
\date{}

\begin{document}

\maketitle

\begin{abstract}
Let $G=(V,E)$ be a hypergraph of rank $r:=\max_{e\in E}|e|.$
A $1/k$-majority $(k+1)$-edge-colouring of $G$ is an assignment of $k+1$ colours to the
edges such that, for every vertex $v$ and every colour $i$, at most
$
\left\lfloor \frac{d(v)}{k}\right\rfloor
$
edges incident with $v$ receive colour $i$.

Motivated by recent work on majority edge-colourings of graphs, we initiate the study of
the corresponding problem for hypergraphs. First, sharpening the probabilistic argument by a $KL$ large-deviation estimate, we obtain
a sufficient minimum-degree condition of order $k^3\log(kr)$ with the sharp large-deviation constant
$
I_k:=D\!\left(\frac1k\middle\|\frac1{k+1}\right)=\Theta(k^{-3}),
$
where $D(\cdot\|\cdot)$ denotes the binary relative entropy.
Our main constructive result shows that
every hypergraph of rank at most $r$ and minimum degree at least $2rk^2$ admits a
$1/k$-majority $(k+1)$-edge-colouring. The proof is based on a hypergraph extension of
the key discrepancy lemma used in the graph case.

We also show that the logarithmic dependence on the rank can be determined asymptotically. If $\mu_k(r)$ denotes the least minimum-degree threshold that guarantees a $1/k$-majority
$(k+1)$-edge-colouring for all hypergraphs of rank at most $r$, then for every fixed $k\ge2$,
$
\mu_k(r)=\frac{\log r}{I_k}+O_k(\log\log r).
$
In particular, the correct logarithmic threshold is of order $k^3\log r$.

Finally, we determine the correct order of the degree--colour trade-off. For integers
$k\ge2$, $p\ge1$, and $r\ge2$, let $\nu_{k,p}(r)$ denote the least integer $q$ such that every
hypergraph of rank at most $r$ and minimum degree at least $kp$ admits a $1/k$-majority
$q$-edge-colouring. Then
$
\nu_{k,p}(r)=\Theta_{k,p}(r^{1/p}).
$
In particular, minimum degree at least $k^2-k$ guarantees a $1/k$-majority
$O_k(r^{1/(k-1)})$-edge-colouring, and this exponent is best possible.
\end{abstract}

\section{Introduction}
Majority colouring is a relaxed form of colouring in which a colour is allowed to appear
locally, but not too often. Its earliest form is the notion of an unfriendly partition,
introduced by Lov\'asz~\cite{lovasz}, which asks for a $2$-colouring of the vertices of an
undirected graph such that each vertex has at most half of its neighbours in its own colour
class. This idea was later studied for infinite graphs~\cite{aharoni}, leading to the well-known
Unfriendly Partition Conjecture.

The concept was subsequently extended to digraphs by Kreutzer, Oum, Seymour,
van der Zypen, and Wood~\cite{kreutzer}, who considered vertex colourings in which no vertex
has more than half of its out-neighbours receiving the same colour. More generally, one may
require that no colour appears on more than a $1/k$ fraction of the relevant neighbourhood,
leading to the notion of a $1/k$-majority colouring. For undirected graphs, the analogous
edge-colouring problem was introduced by Bock et al.~\cite{bock}. They proved that every graph
with minimum degree at least $2$ admits a majority $4$-edge-colouring, and more generally
that graphs of sufficiently large minimum degree admit $1/k$-majority $(k+1)$-edge-colourings.

P\k{e}ka{\l}a and Przyby{\l}o~\cite{pekala} substantially improved the graph case by proving that
minimum degree $2k^2$ is sufficient, and they conjectured that the optimal threshold is $k^2$.
Their argument is based on a discrepancy lemma that rounds edge weights while preserving
strong local control. This naturally raises the question of how the problem behaves for
hypergraphs.

In this paper, we initiate the study of majority edge-colourings in hypergraphs.
Let $G=(V,E)$ be a hypergraph, and let
$
r:=\max_{e\in E}|e|
$
be its rank. A $1/k$-majority $(k+1)$-edge-colouring of $G$ is an assignment of
$k+1$ colours to the edges such that, for every vertex $v$ and every colour $i$, at most
$
\left\lfloor \frac{d(v)}{k}\right\rfloor
$
edges incident with $v$ receive colour $i$.

Our first result is a probabilistic existence theorem. Using the weighted Lov\'asz Local
Lemma together with the exact KL large-deviation rate, we obtain a sufficient minimum-degree
condition of order $k^3\log(kr)$ for a $1/k$-majority $(k+1)$-edge-colouring.
Our main constructive result shows that every hypergraph of rank at most $r$ and minimum
degree at least $2rk^2$ admits such a colouring.

We then turn to the threshold function $\mu_k(r)$. A self-contained random-family construction,
proved in Appendix~\ref{app:mu-lower}, shows that for every fixed integer $k\ge2$,
$$
\mu_k(r)=\frac{\log r}{I_k}+O_k(\log\log r),
\qquad
I_k:=D\!\left(\frac1k\middle\|\frac1{k+1}\right),
$$
where
$$
D(a\|b):=a\log\frac{a}{b}+(1-a)\log\frac{1-a}{1-b}
\qquad (0<a,b<1)
$$
denotes the binary relative entropy. Thus the logarithmic dependence on the rank is
asymptotically sharp, including the leading constant.

Finally, we study the effect of allowing more colours. If $\nu_{k,p}(r)$ denotes the least
number of colours that always suffices for rank-$r$ hypergraphs of minimum degree at least
$kp$, then we prove $
\nu_{k,p}(r)=\Theta_{k,p}(r^{1/p}).$
This determines the correct order of the degree--colour trade-off. In particular,
minimum degree at least $k^2-k$ already guarantees a $1/k$-majority
$O_k(r^{1/(k-1)})$-edge-colouring.

For clarity, we summarise our main contributions below.

\begin{itemize}[leftmargin=2em]
\item We obtain, via the weighted Lov\'asz Local Lemma and a KL large-deviation estimate,
  a sufficient minimum-degree threshold of order $k^3\log(kr)$ for the existence of a
  $1/k$-majority $(k+1)$-edge-colouring.

\item We prove that every hypergraph $G$ of rank at most $r$ and minimum degree $\delta(G)\ge 2rk^2$ admits a $1/k$-majority $(k+1)$-edge-colouring.

\item We determine the asymptotic logarithmic threshold: for every fixed integer $k\ge2$,
  $$
  \mu_k(r)=\frac{\log r}{I_k}+O_k(\log\log r),
  \qquad
  I_k=D\!\left(\frac1k\middle\|\frac1{k+1}\right).
  $$

\item We obtain a general $q$-colour obstruction: for all integers $k\ge2$, $q\ge1$,
  and $d\ge k$, there exists a $d$-regular hypergraph of rank
  $
  \binom{q\lfloor d/k\rfloor}{\lfloor d/k\rfloor}
  $
  with no $1/k$-majority $q$-edge-colouring.

\item We determine the correct fixed-$p$ degree--colour trade-off:
  $$
  \nu_{k,p}(r)=\Theta_{k,p}(r^{1/p})
  \qquad (k\ge2,\ p\ge1).
  $$
  In particular, minimum degree at least $k^2-k$ guarantees a $1/k$-majority
  $O_k(r^{1/(k-1)})$-edge-colouring.
\end{itemize}

The remainder of the paper is organised as follows. In Section~\ref{sec:preliminaries}, we introduce the basic
definitions and notation. In Section~\ref{sec:main-results}, we prove the probabilistic threshold,
establish the hypergraph discrepancy lemma, derive the main theorem, present the counterexample
showing that rank dependence is unavoidable, and obtain a general degree--colour trade-off.
In Section~\ref{sec:threshold-function}, we introduce the threshold function $\mu_k(r)$ and show
the fixed-$k$ asymptotic consequence of our results. Appendix~\ref{app:mu-lower} contains the
self-contained random-family argument proving the matching lower bound for $\mu_k(r)$.
We conclude with some further remarks and open problems.

\section{Preliminaries}\label{sec:preliminaries}

\begin{definition}[Hypergraph]
A hypergraph $G=(V,E)$ consists of a finite set $V$ of vertices and a collection $E$ of
non-empty subsets of $V$, called hyperedges. The degree of a vertex $v\in V$, denoted by
$d_G(v)$, is the number of hyperedges containing $v$. The minimum degree is
$$
\delta(G):=\min_{v\in V} d_G(v).
$$
The quantity
$$
\rank(G):=\max_{e\in E}|e|
$$
is called the rank of $G$.
\end{definition}

\begin{definition}[$1/k$-majority edge-colouring]
Let $k\ge 2$ and $q\ge 1$ be integers. A $1/k$-majority $q$-edge-colouring of a hypergraph
$G=(V,E)$ is a function
$$
c:E\to \{1,2,\dots,q\}
$$
such that for every vertex $v\in V$ and every colour $i\in\{1,2,\dots,q\}$,
$$
\bigl|\{e\in E: v\in e \text{ and } c(e)=i\}\bigr|\le \left\lfloor \frac{d_G(v)}{k}\right\rfloor.
$$
\end{definition}

For $0<a,b<1$, we write
$$
D(a\|b):=a\log\frac{a}{b}+(1-a)\log\frac{1-a}{1-b}
$$
for the binary relative entropy.
\begin{definition}
For integers $k,r\ge 2$, let $\mu_k(r)$ denote the least integer $d$ such that every
hypergraph of rank at most $r$ and minimum degree at least $d$ admits a
$1/k$-majority $(k+1)$-edge-colouring.
\end{definition}
\section{Main results}\label{sec:main-results}

\subsection{A probabilistic threshold}

Under essentially the same framework as Bock et al.~\cite[Theorem~4]{bock}, one can obtain
an existence threshold for hypergraphs. For completeness, we include the full proof.

\begin{theorem}\label{thm:prob-threshold}
For all integers $k\ge 2$ and $r\ge 2$, every hypergraph of rank at most $r$ and
minimum degree at least $\delta_{k,r}$ admits a $1/k$-majority $(k+1)$-edge-colouring,
provided that $\delta_{k,r}$ satisfies the inequalities
$$
4(k+1)e^{-\delta_{k,r} I_k}\le 1,
$$
and
$$
8(k+1)(r-1)e^{-\delta_{k,r} I_k}\,\delta_{k,r}\le 1,
$$ where $$I_k:=D\!\left(\frac1k\middle\|\frac1{k+1}\right).$$
\end{theorem}

\begin{proof}
Let $G$ be a hypergraph of minimum degree $\delta\ge \delta_{k,r}$, and let
$$
c:E(G)\to [k+1]
$$
be a random $(k+1)$-edge-colouring obtained by choosing the colour of each edge uniformly
and independently at random.

For every vertex $u\in V(G)$, let $A_u$ be the bad event that more than $d_G(u)/k$ of the
edges incident with $u$ receive the same colour. Writing $d=d_G(u)$, the union bound and the
Chernoff inequality (see, for instance,~\cite{molloy}) gives
\begin{align*}
\mathbb{P}[A_u]
&\le (k+1)\,\mathbb{P}\!\left[\operatorname{Bin}\!\left(d,\frac1{k+1}\right)>\frac{d}{k}\right]\\
&=(k+1)\,\mathbb{P}\!\left[\operatorname{Bin}\!\left(d,\frac1{k+1}\right)>
\left(1+\frac1k\right)\frac{d}{k+1}\right]\\
&\le (k+1)e^{-d I_k}.
\end{align*}

Set
$$
p:=(k+1)e^{-\delta I_k},
\qquad
t_u:=\left\lfloor \frac{d_G(u)}{\delta}\right\rfloor.
$$
Since $d_G(u)\ge \delta$, the integer $t_u$ is positive, and
$$
2t_u = 2\left\lfloor \frac{d_G(u)}{\delta}\right\rfloor \ge \frac{d_G(u)}{\delta}.
$$
By the first displayed inequality, we have $p\le 1/4$, and hence
$$
\mathbb{P}[A_u]\le p^{d_G(u)/\delta}\le p^{t_u}.
$$

Let $N_G(u)$ denote the set of vertices $v\neq u$ for which there exists an edge containing both
$u$ and $v$. The event $A_u$ depends only on the colours of the edges incident with $u$, and
is therefore mutually independent of all events $A_v$ with $v\notin N_G(u)$. Moreover,
$$
|N_G(u)|\le (r-1)d_G(u),
$$
since each edge incident with $u$ contains at most $r-1$ other vertices. Because $t_v\ge 1$ for
every $v$, we obtain
\begin{align*}
\sum_{v\in N_G(u)} (2p)^{t_v}
&\le 2p\,|N_G(u)|\\
&\le 2p(r-1)d_G(u)\\
&\le 4p\delta(r-1)t_u\\
&=\Bigl(4(k+1)(r-1)e^{-\delta I_k}\,\delta\Bigr)t_u,
\end{align*}
which is at most $t_u/2$ by the second displayed inequality.

Lemma~\ref{lem:wlll}, proved below, now implies that with positive probability none of the
events $A_u$ occurs. Therefore $G$ admits a $1/k$-majority $(k+1)$-edge-colouring.
\end{proof}

\subsection{A hypergraph discrepancy lemma}

We next recall the key graph-theoretic rounding lemma of P\k{e}ka{\l}a and
Przyby{\l}o~\cite[Lemma~6]{pekala}.

\begin{lemma}[{\cite[Lemma~6]{pekala}}]\label{lem:graph-discrepancy}
Let $G=(V,E)$ be a graph, and let $z:E\to [0,1]$ be a weight function. Then there exists a
function $x:E\to\{0,1\}$ such that:
\begin{enumerate}[label=\textnormal{(\roman*)},leftmargin=2em]
\item
$$
\sum_{e\ni v} z(e)-1 < \sum_{e\ni v} x(e) \le \sum_{e\ni v} z(e)+1
\qquad\text{for every }v\in V;
$$
\item if
$$
\sum_{e\ni u} x(e)<\sum_{e\ni u} z(e)
\quad\text{and}\quad
\sum_{e\ni v} x(e)<\sum_{e\ni v} z(e)
$$
for some edge $uv\in E$, then $x(uv)=1$;
\item each vertex $v$ with
$$
\sum_{e\ni v} x(e)=\sum_{e\ni v} z(e)+1
$$
belongs to an odd cycle $C_v$ such that
$$
\sum_{e\ni u} z(e)\in \mathbb{Z}
\qquad\text{for every }u\in V(C_v),
$$
and moreover all such cycles are pairwise vertex-disjoint.
\end{enumerate}
\end{lemma}

We extend the first part of Lemma~\ref{lem:graph-discrepancy} to hypergraphs. Our proof
follows the same general constructive philosophy, but the larger edge sizes force an additive
error depending on the rank.

\begin{lemma}\label{lem:hypergraph-discrepancy}
Let $G=(V,E)$ be a hypergraph of rank at most $r$, and let $z:E\to [0,1]$ be a function.
Then there exists a function $x:E\to\{0,1\}$ such that
$$
\sum_{e\ni v} z(e)-r < \sum_{e\ni v} x(e) < \sum_{e\ni v} z(e)+r
\qquad\text{for every }v\in V.
$$
\end{lemma}

\begin{proof}
For every edge $e$ with $z(e)\in\{0,1\}$, simply set $x(e)=z(e)$ and delete $e$. Thus it
suffices to consider the remaining hypergraph $H$, in which every edge $e$ satisfies
$0<z(e)<1$.

If every vertex of $H$ has degree at most $r$, then the claim is immediate: for every vertex
$v$, the sum \(\sum_{e\ni v} x(e)\) involves at most $r$ terms, each lying in $\{0,1\}$,
whereas \(\sum_{e\ni v} z(e)\) involves the same number of terms, each lying in $(0,1)$.
Hence
$$
\left|\sum_{e\ni v} x(e)-\sum_{e\ni v} z(e)\right|<r
$$
for every choice of $x:E(H)\to\{0,1\}$.

We may therefore assume that some vertices have degree at least $r+1$. Let
$$
U=\{v_1,\dots,v_s\}
$$
be the set of vertices of $H$ whose degree is at least $r+1$, and consider the polytope
$$
\mathcal{P}=
\left\{h\in [0,1]^{E(H)}:
\sum_{e\ni v_i} h(e)=\sum_{e\ni v_i} z(e)
\text{ for every }i\in[s]
\right\}.
$$
Since $z\in \mathcal{P}$, the polytope is nonempty. Choose an extreme point $h$ of
$\mathcal{P}$.

We claim that some coordinate of $h$ belongs to $\{0,1\}$. Indeed, suppose to the contrary
that every coordinate of $h$ lies strictly between $0$ and $1$. Then the only active
constraints at $h$ are the $s$ linear equations defining $\mathcal{P}$. On the other hand, the
number of variables is $|E(H)|$, and
$$
(r+1)s \le \sum_{i=1}^s d_H(v_i)
\le \sum_{e\in E(H)} |e| \le r|E(H)|.
$$
Hence $s<|E(H)|$. Therefore the homogeneous system
$$
\sum_{e\ni v_i} y(e)=0
\qquad\text{for every }i\in[s]
$$
has a nonzero solution $y\in \mathbb{R}^{E(H)}$. For sufficiently small $\varepsilon>0$, both
$h+\varepsilon y$ and $h-\varepsilon y$ still belong to $[0,1]^{E(H)}$, and they also satisfy all
defining equalities of $\mathcal{P}$. This contradicts the extremality of $h$.

Thus some edge $e_1$ satisfies $h(e_1)\in\{0,1\}$. Set $x(e_1)=h(e_1)$, delete $e_1$, and
replace $z$ by $h$ on the remaining edges. For every vertex in $U$, the total weight of the
remaining incident edges is unchanged, so the same argument may be repeated. Since at least
one edge is removed at each step, the process terminates after finitely many steps. At the end,
every remaining vertex has degree at most $r$, and we may assign the remaining edges
arbitrarily.

For a vertex $v\in V(H)$, let $T(v)$ be the first step at which the current residual degree of
$v$ becomes at most $r$; if $v$ already has degree at most $r$ initially, set $T(v)=0$. Up to
step $T(v)-1$, the vertex $v$ belongs to the constrained set $U$, so at each of those steps the
sum of the already fixed values on deleted edges incident with $v$ together with the current
residual weights on the remaining incident edges is preserved exactly. At time $T(v)$, the
residual hypergraph contains at most $r$ edges through $v$. From that point onward we no longer
need exact preservation: every further change at $v$ comes only from these residual edges, and
each such edge contributes either a number in $(0,1)$ to the residual weighted sum at time
$T(v)$ or a value in $\{0,1\}$ to the final rounding. Consequently, the total discrepancy created
after time $T(v)$ is strictly smaller than $r$. Therefore
$$
\left|\sum_{e\ni v} x(e)-\sum_{e\ni v} z(e)\right|<r
$$
for every vertex $v$, as required.
\end{proof}

\subsection{Proof of the main theorem}

The proof of our main theorem adapts the iterative partitioning framework from
P\k{e}ka{\l}a and Przyby{\l}o~\cite[Theorem~12]{pekala}. The only substantial change is the
choice of the parameter $\alpha_i$ in each round, which must absorb the larger discrepancy
from Lemma~\ref{lem:hypergraph-discrepancy}.

\begin{theorem}\label{thm:main}
Let $G=(V,E)$ be a hypergraph of rank at most $r$. If
$$
\delta(G)\ge 2rk^2,
$$
then $G$ admits a $1/k$-majority $(k+1)$-edge-colouring.
\end{theorem}

\begin{proof}
Set $\delta:=\delta(G)$ and define $H_0:=G$. For $i=1,\dots,k$, let
$$
S_i:=\delta-i\left(\frac{\delta}{k}-2r\right),
\qquad
\alpha_i:=\frac{\delta/k-r}{S_{i-1}}.
$$
Since
$$
S_{i-1}\ge \frac{\delta}{k}+2r(k-i+1)>\frac{\delta}{k}-r>0,
$$
we have $0<\alpha_i<1$ for every $i\in[k]$. Apply Lemma~\ref{lem:hypergraph-discrepancy} to
$H_{i-1}$ with the constant weight function $z\equiv \alpha_i$. This yields a spanning subhypergraph $C_i\subseteq H_{i-1}$ such that
$$
\alpha_i d_{H_{i-1}}(v)-r < d_{C_i}(v) < \alpha_i d_{H_{i-1}}(v)+r
\qquad\text{for every }v\in V.
$$
Now define
$$
H_i:=H_{i-1}-E(C_i).
$$
Finally, let $C_{k+1}:=H_k$. Then $C_1,\dots,C_{k+1}$ form an edge-partition of $G$.
We show that each colour class satisfies the $1/k$-majority condition.

For each vertex $v\in V$, write
$$
\beta_v:=\frac{d_G(v)}{\delta}\ge 1.
$$
We claim that for every $i\in\{0,1,\dots,k\}$,
$$
d_{H_i}(v)\le \beta_v S_i
\qquad\text{for every }v\in V.
$$
The case $i=0$ is immediate, since $S_0=\delta$.
Assume the claim holds for $i-1$. Then
\begin{align*}
d_{H_i}(v)
&= d_{H_{i-1}}(v)-d_{C_i}(v)\\
&< (1-\alpha_i)d_{H_{i-1}}(v)+r\\
&\le (1-\alpha_i)\beta_v S_{i-1}+r\\
&= \beta_v\Bigl(S_{i-1}-\bigl(\tfrac{\delta}{k}-r\bigr)\Bigr)+r\\
&\le \beta_v\Bigl(S_{i-1}-\tfrac{\delta}{k}+2r\Bigr)\\
&= \beta_v S_i,
\end{align*}
where we used $\alpha_iS_{i-1}=\delta/k-r$ and $r\le \beta_v r$ in the penultimate line.
This proves the claim.

We now estimate the colour classes $C_i$ for $1\le i\le k$:
\begin{align*}
d_{C_i}(v)
&< \alpha_i d_{H_{i-1}}(v)+r\\
&\le \alpha_i\beta_v S_{i-1}+r\\
&= \beta_v\Bigl(\frac{\delta}{k}-r\Bigr)+r\\
&\le \beta_v\frac{\delta}{k}
 = \frac{d_G(v)}{k}.
\end{align*}
Thus every one of the first $k$ colour classes satisfies the $1/k$-majority condition.

Finally,
$$
d_{C_{k+1}}(v)=d_{H_k}(v)\le \beta_v S_k = \beta_v\cdot 2rk.
$$
Since $\delta\ge 2rk^2$, we have $2rk\le \delta/k$, and therefore
$$
d_{C_{k+1}}(v)\le \beta_v\frac{\delta}{k}=\frac{d_G(v)}{k}.
$$
Hence $C_1,\dots,C_{k+1}$ form a $1/k$-majority $(k+1)$-edge-colouring of $G$.
\end{proof}

\subsection{Necessity of rank dependence}

\begin{remark}\label{rem:rank-necessary}
The dependence on the rank in our bounds is necessary. In particular, there is no function of
$k$ alone that guarantees the existence of a $1/k$-majority $(k+1)$-edge-colouring for all
hypergraphs.
\end{remark}

\begin{theorem}\label{thm:regular-counterexample-general-q}
Let $k\ge 2$, $q\ge 1$, and $d\ge k$ be integers, and let
$
m:=\left\lfloor \frac{d}{k}\right\rfloor.
$
Then there exists a $d$-regular hypergraph $H_{k,q,d}$ such that
$
\rank(H_{k,q,d})=\binom{qm}{m},
$
and $H_{k,q,d}$ does not admit a $1/k$-majority $q$-edge-colouring.
\end{theorem}

\begin{proof}
Let $N:=qm+1$, and take $N$ shared edges $e_1,\dots,e_N$.
For every $(m+1)$-subset $X\subseteq [N]$, create a vertex $v_X$ incident with precisely the
shared edges $\{e_i:i\in X\}$ together with $d-(m+1)$ private edges that meet no other vertex.
This is possible because $m+1\le d$ for $k\ge2$.

Every vertex has degree $d$, so the resulting hypergraph is $d$-regular. Every private edge has
size $1$, while a shared edge $e_i$ is incident with exactly those vertices $v_X$ for which
$i\in X$, and therefore
$$
|e_i|=\binom{N-1}{m}=\binom{qm}{m}.
$$
Hence $\rank(H_{k,q,d})=\binom{qm}{m}$.

Now consider any $q$-colouring of the shared edges. Since there are $qm+1$ shared edges and
only $q$ colours, some colour appears on at least $m+1$ shared edges. Let $X$ be an
$(m+1)$-subset of indices of shared edges receiving that colour. The corresponding vertex $v_X$
is incident with $m+1$ edges of one colour, but
$$
m+1>\left\lfloor \frac{d}{k}\right\rfloor,
$$
so the $1/k$-majority condition fails at $v_X$.
\end{proof}

\subsection{More colours and a degree--colour trade-off}

The next result shows that allowing more colours substantially lowers the minimum-degree
requirement, and that no sparsity assumption is needed.

\begin{theorem}\label{thm:many-colours-general}
Every hypergraph of rank at most $r$ and minimum degree at least $k^2-k$
admits a $1/k$-majority $(kr+1)$-edge-colouring.
\end{theorem}

\begin{proof}
Let $G$ be a hypergraph of rank at most $r$ and minimum degree at least $k^2-k$.
For each vertex $u\in V(G)$, write
$$
d_G(u)=t_uk+m_u,
\qquad 0\le m_u\le k-1,
$$
so that
$$
t_u=\left\lfloor \frac{d_G(u)}{k}\right\rfloor.
$$
Since $d_G(u)\ge k^2-k$, we have $t_u\ge k-1$, and hence $m_u\le t_u$.

We split each vertex $u$ into $t_u$ new vertices
$$
u_1,\dots,u_{t_u},
$$
and partition the edges incident with $u$ into $t_u$ groups in such a way that
exactly $m_u$ groups have size $k+1$ and the remaining $t_u-m_u$ groups have size $k$.
For every edge $e=\{u^{(1)},\dots,u^{(s)}\}\in E(G)$, if the incidence of $e$ at $u^{(j)}$
is assigned to the copy $u^{(j)}_{i_j}$, then we replace $e$ by the edge
$$
e^*=\{u^{(1)}_{i_1},\dots,u^{(s)}_{i_s}\}.
$$
In this way we obtain a hypergraph $G^*$ whose edge set is naturally identified with $E(G)$.

By construction,
$$
\Delta(G^*)\le k+1,
$$
and every edge of $G^*$ still has size at most $r$.

Now consider the line graph $L(G^*)$. Its vertices are the edges of $G^*$, and two
vertices of $L(G^*)$ are adjacent whenever the corresponding edges of $G^*$ intersect.
Let $e^*\in E(G^*)$. Since $e^*$ contains at most $r$ vertices and every vertex of $G^*$
has degree at most $k+1$, the edge $e^*$ can meet at most
$$
r((k+1)-1)=kr
$$
other edges. Therefore,
$$
\Delta(L(G^*))\le kr.
$$
Hence $L(G^*)$ admits a proper vertex-colouring with at most $kr+1$ colours by a greedy
algorithm. Equivalently, $G^*$ admits a proper edge-colouring with at most $kr+1$ colours.

Pull this colouring back to the original hypergraph $G$. Fix a vertex $u\in V(G)$.
A given colour appears at most once on the edges incident with each split copy
$u_1,\dots,u_{t_u}$, because the colouring on $G^*$ is proper. Therefore, in $G$ the same
colour appears at most $t_u$ times on the edges incident with $u$. Since
$$
t_u=\left\lfloor \frac{d_G(u)}{k}\right\rfloor,
$$
the colouring is a $1/k$-majority $(kr+1)$-edge-colouring of $G$.
\end{proof}

The same splitting argument yields a more general degree--colour trade-off.

\begin{definition}
For integers $k\ge2$, $p\ge1$, and $r\ge2$, let $\nu_{k,p}(r)$ denote the least integer $q$
such that every hypergraph of rank at most $r$ and minimum degree at least $kp$ admits a
$1/k$-majority $q$-edge-colouring.
\end{definition}
\begin{lemma}\label{lem:wlll}
Let $\{A_u:u\in U\}$ be events. Suppose that for each $u\in U$ we are given:
\begin{itemize}[leftmargin=1.6em]
\item a positive integer $t_u$,
\item a set $N(u)\subseteq U\setminus\{u\}$ such that $A_u$ is mutually independent of all events $A_v$ with $v\notin N(u)$,
\item a real number $x\in(0,1/2]$ such that
$$
\Pr(A_u)\le x^{t_u}
\quad\text{and}\quad
\sum_{v\in N(u)} (2x)^{t_v}\le \frac{t_u}{2}
$$
for every $u\in U$.
\end{itemize}
Then $\Pr\bigl(\bigcap_{u\in U}\overline{A_u}\bigr)>0$.
\end{lemma}

\begin{proof}
Set
$$
y_u:=(2x)^{t_u}
\qquad (u\in U).
$$
Since $x\le 1/2$ and $t_u\ge1$, we have $0<y_u\le 1/2$ for every $u$. Moreover,
$$
\Pr(A_u)\le x^{t_u}=2^{-t_u}y_u.
$$
For $0\le y\le 1/2$ we have $1-y\ge 2^{-2y}$, and therefore
$$
\prod_{v\in N(u)}(1-y_v)\ge 2^{-2\sum_{v\in N(u)}y_v}\ge 2^{-t_u}.
$$
Consequently,
$$
\Pr(A_u)\le y_u\prod_{v\in N(u)}(1-y_v)
\qquad\text{for every }u\in U.
$$
The asymmetric Lov\'asz Local Lemma now implies that
$\Pr\bigl(\bigcap_{u\in U}\overline{A_u}\bigr)>0$.
\end{proof}
\begin{theorem}\label{thm:degree-colour-tradeoff}
Let $k\ge2$, $p\ge1$, and $r\ge2$ be integers. Every hypergraph of rank at most $r$ and
minimum degree at least $kp$ admits a $1/k$-majority $q$-edge-colouring whenever
$$
q^p\ge 8e^{p+1}k^{p+2}p\,(r-1).
$$
In particular,
$$
\nu_{k,p}(r)\le \left\lceil\bigl(8e^{p+1}k^{p+2}p\,(r-1)\bigr)^{1/p}\right\rceil.
$$
\end{theorem}
\begin{proof}
Let $G$ be a hypergraph of rank at most $r$ and minimum degree at least $kp$. Colour each edge independently and uniformly from $[q]:=\{1,\dots,q\}$.

For a vertex $u\in V(G)$, write $d=d_G(u)$ and define
$$
s_u:=\floor{\frac{d}{kp}},
\qquad
\tau_u:=\floor{\frac{d}{k}}+1.
$$
Since $d\ge kp$, the integer $s_u$ is positive. Let $A_u$ be the bad event that some colour appears on at least $\tau_u$ edges incident with $u$.

Fix a colour $i\in[q]$. The number of incident edges of colour $i$ at $u$ has distribution $\Bin(d,1/q)$. Using the inequality
$$
\mathbf 1_{\{X\ge t\}}\le \binom{X}{t}
\qquad (X\in\mathbb N,\ t\in\mathbb N),
$$
we get
$$
\Pr\!\left[\Bin\!\left(d,\frac1q\right)\ge \tau_u\right]
\le
\mathbb E\!\left[\binom{\Bin(d,1/q)}{\tau_u}\right]
=
\binom{d}{\tau_u}q^{-\tau_u}.
$$
Hence, by the union bound,
$$
\Pr(A_u)\le q\binom{d}{\tau_u}q^{-\tau_u}.
$$
Now $\tau_u>d/k$, so $d/\tau_u<k$, and therefore
$$
\binom{d}{\tau_u}\le \left(\frac{ed}{\tau_u}\right)^{\tau_u}< (ek)^{\tau_u}.
$$
Also,
$$
\tau_u=\floor{\frac{d}{k}}+1\ge p\floor{\frac{d}{kp}}+1 = ps_u+1.
$$
Thus
$$
\Pr(A_u)\le q\left(\frac{ek}{q}\right)^{\tau_u}.
$$
Our hypothesis implies
$$
q^p\ge 8e^{p+1}k^{p+2}p(r-1)>(ek)^p,
$$
and hence $q>ek$. Therefore $ek/q<1$, and since $\tau_u\ge ps_u+1$ we obtain
$$
\Pr(A_u)
\le q\left(\frac{ek}{q}\right)^{ps_u+1}.
$$
Set
$$
x:=\frac{(ek)^{p+1}}{q^p}.
$$
Since $s_u\ge1$, the last inequality implies $\Pr(A_u)\le x^{s_u}$.

Next define $N(u)$ to be the set of vertices $v\neq u$ for which some edge contains both $u$ and $v$. Then $A_u$ depends only on the colours of the edges incident with $u$, so $A_u$ is mutually independent of all $A_v$ with $v\notin N(u)$. Moreover,
$$
|N(u)|\le (r-1)d,
$$
since every edge through $u$ contains at most $r-1$ further vertices. Because $d\ge kp$, we also have
$$
2s_u = 2\floor{\frac{d}{kp}}\ge \frac{d}{kp},
\qquad\text{so}\qquad d\le 2kp\,s_u.
$$
Therefore, provided $x\le 1/4$,
$$
\sum_{v\in N(u)} (2x)^{s_v}
\le 2x|N(u)|
\le 2x(r-1)d
\le 4x(r-1)kp\,s_u.
$$
Our hypothesis on $q$ is exactly
$$
8x(r-1)kp\le 1.
$$
Since $r\ge2$ and $kp\ge2$, this also implies $x\le1/4$. Hence
$$
\sum_{v\in N(u)} (2x)^{s_v}\le \frac{s_u}{2}.
$$
Lemma~\ref{lem:wlll} now shows that with positive probability none of the events $A_u$ occur.
\end{proof}

\begin{corollary}\label{cor:tradeoff-lower}
For all fixed integers $k\ge2$ and $p\ge1$,
$$
\left\lfloor \frac{(p!\,r)^{1/p}}{p}\right\rfloor+1
\le
\nu_{k,p}(r)
\le
\left\lceil\bigl(8e^{p+1}k^{p+2}p\,(r-1)\bigr)^{1/p}\right\rceil.
$$
In particular,
$$
\nu_{k,p}(r)=\Theta_{k,p}(r^{1/p}).
$$
\end{corollary}

\begin{proof}
The upper bound is Theorem~\ref{thm:degree-colour-tradeoff}.
For the lower bound, apply Theorem~\ref{thm:regular-counterexample-general-q} with $d=kp$.
Then $m=p$, so for every integer $q\ge1$ there exists a $kp$-regular hypergraph of rank
$
\binom{qp}{p}
$
with no $1/k$-majority $q$-edge-colouring. Hence, whenever $\binom{qp}{p}\le r$, one must have
$\nu_{k,p}(r)\ge q+1$. Since
$
\binom{qp}{p}\le \frac{(qp)^p}{p!},
$
it is enough to choose
$
q=\left\lfloor \frac{(p!\,r)^{1/p}}{p}\right\rfloor.
$
\end{proof}

\begin{corollary}
If the minimum degree is at least $k^2-k$, then every hypergraph of rank at most $r$ admits a
$1/k$-majority $O_k(r^{1/(k-1)})$-edge-colouring. Moreover, the exponent $1/(k-1)$ is best
possible.
\end{corollary}

While the probabilistic argument needs the prerequisite $q > ek$, we can extend Theorem~\ref{thm:main} and get a bound when arbitrary number of additional colours are permitted.
\begin{theorem} \label{thm:generalized_peeling}
Let $k \ge 2$, $q > k$, and $r \ge 2$ be integers. Every hypergraph $G$ of rank at most $r$ and minimum degree
$$ \delta(G) \ge \frac{2rk(q-1)}{q-k} $$
admits a $1/k$-majority $q$-edge-colouring.
\end{theorem}

\begin{proof}
Let $\delta := \delta(G)$. We adopt the iterative peeling procedure and notations established in the proof of Theorem ~\ref{thm:main}. Specifically, for each vertex $v \in V(G)$, let $\beta_v := d_G(v)/\delta \ge 1$. We define the sequence of residual degree bounds $S_0 := \delta$ and $S_i := S_{i-1} - (\frac{\delta}{k} - 2r)$ for $i \ge 1$. 

At each step $i$ (for $i = 1, \dots, q-1$), we evaluate the current degree bound $S_{i-1}$ of the residual hypergraph $H_{i-1}$:

\textbf{Case 1 (Early Stopping):} If $S_{i-1} \le \frac{\delta}{k}$, we simply assign all remaining edges to the current class, $C_i := E(H_{i-1})$, and leave all subsequent colour classes empty ($C_j := \emptyset$ for $j > i$). Because $d_{C_i}(v) \le \beta_v S_{i-1} \le \beta_v \frac{\delta}{k} = \frac{d_G(v)}{k}$, the integer degree trivially satisfies $d_{C_i}(v) \le \lfloor d_G(v)/k \rfloor$. The empty classes also trivially satisfy the $1/k$-majority condition, and the algorithm terminates successfully.

\textbf{Case 2 (Greedy Peeling):} If $S_{i-1} > \frac{\delta}{k}$, we extract $C_i$ from $H_{i-1}$ via Lemma ~\ref{lem:hypergraph-discrepancy} exactly as in Theorem ~\ref{thm:main}. The required fractional weight $\alpha_i := (\delta/k - r)/S_{i-1}$ is strictly positive because our hypothesis $\delta \ge \frac{2rk(q-1)}{q-k} > kr$ guarantees $\delta/k > r$. As established in Theorem ~\ref{thm:main}, this yields a valid colour class $C_i$ satisfying $d_{C_i}(v) \le \lfloor d_G(v)/k \rfloor$, and the updated residual hypergraph $H_i$ satisfies $d_{H_i}(v) \le \beta_v S_i$.

It remains to verify that if the algorithm performs $q-1$ full peeling steps, the final residual hypergraph $H_{q-1}$ can act as the last valid colour class $C_q$. This requires $S_{q-1} \le \frac{\delta}{k}$. By unrolling the recurrence, we have $S_{q-1} = \delta - (q-1)(\frac{\delta}{k} - 2r)$. Imposing $S_{q-1} \le \frac{\delta}{k}$ gives:
$$ \delta - \frac{\delta}{k} + 2r(q-1) \le \frac{q-1}{k}\delta. $$
Rearranging this inequality yields exactly $\delta \left( \frac{q-k}{k} \right) \ge 2r(q-1)$, which matches our initial hypothesis. Thus, the remaining edges can always be legally assigned to $C_q$, completing the constructive proof.
\end{proof}

\section{The threshold function \texorpdfstring{$\mu_k(r)$}{mu\_k(r)}}\label{sec:threshold-function}

We first recall the definition of $\mu_k(r)$.
\begin{definition}
For integers $k,r\ge 2$, let $\mu_k(r)$ denote the least integer $d$ such that every
hypergraph of rank at most $r$ and minimum degree at least $d$ admits a
$1/k$-majority $(k+1)$-edge-colouring.
\end{definition}

\begin{theorem}\label{thm:mu-asymptotic}
For every fixed integer $k\ge2$,
$$
\mu_k(r)=\frac{\log r}{I_k}+O_k(\log\log r)
\qquad (r\to\infty),
$$
where
$
I_k:=D\!\left(\frac1k\middle\|\frac1{k+1}\right).
$
\end{theorem}

\begin{proof}
The upper bound follows from Theorem~\ref{thm:prob-threshold} via the inequality
$
\mu_k(r)\le \frac{\log r}{I_k}+O_k(\log\log r).
$
For the lower bound, Appendix~\ref{app:mu-lower} proves that for all sufficiently large
multiples $d$ of $k$, there exists a $d$-regular hypergraph with no $1/k$-majority
$(k+1)$-edge-colouring and rank at most
$
\exp\bigl(dI_k+C_k\log d\bigr)
$
for a constant $C_k$ depending only on $k$. Let $d$ be the largest multiple of $k$ satisfying
$$
dI_k+C_k\log d\le \log r.
$$
Then $d=\frac{\log r}{I_k}-O_k(\log\log r)$, and the cited lower bound gives
$
\mu_k(r)\ge d\ge \frac{\log r}{I_k}-O_k(\log\log r).
$
Combining the upper and lower bounds yields the theorem.
\end{proof}

\begin{remark}
Since
$$
I_k=\frac1{2k^3}-\frac1{6k^4}+O\!\left(\frac1{k^5}\right),
$$
Theorem~\ref{thm:mu-asymptotic} shows that the logarithmic threshold has leading scale
$2k^3\log r$.
\end{remark}
\section*{Concluding remarks}

For integers $k,r\ge2$, let $\mu_k(r)$ denote the least integer $d$ such that every hypergraph
of rank at most $r$ and minimum degree at least $d$ admits a $1/k$-majority $(k+1)$-edge-colouring.
By Theorem~\ref{thm:mu-asymptotic}, for every fixed integer $k\ge2$,
$$
\mu_k(r)=\frac{\log r}{D(1/k\,\|\,1/(k+1))}+O_k(\log\log r).
$$
Thus the logarithmic rank threshold is now asymptotically determined, including the sharp
large-deviation constant.

On the other hand, the constructive theorem
$
\delta(G)\ge 2rk^2
$
remains the main unresolved $(k+1)$-colour bound. It would be especially interesting to know
whether the factor $2$ can be removed, or more generally whether one can obtain a constructive
upper bound closer to the logarithmic threshold in sparse rank regimes.

Finally, Corollary~\ref{cor:tradeoff-lower} determines the correct order of the fixed-$p$
degree--colour trade-off:
$$
\nu_{k,p}(r)=\Theta_{k,p}(r^{1/p}).
$$
So the main remaining questions now concern sharp constants, constructive versions of the
polynomial trade-off, and possible refinements under additional structural assumptions such as
small codegree or intersection sparsity.

\appendix

\section{A self-contained lower bound for \texorpdfstring{$\mu_k(r)$}{mu\_k(r)}}\label{app:mu-lower}

In this appendix we prove the lower bound used in Theorem~\ref{thm:mu-asymptotic}. Fix an
integer $k\ge2$, let
$$
q:=k+1,
\qquad
I_k:=D\!\left(\frac1k\middle\|\frac1{k+1}\right),
$$
and write $d=km$ with $m\in\mathbb N$.

We first use the standard incidence duality between uniform set systems and regular
hypergraphs.

\begin{lemma}\label{lem:incidence-duality}
Let $\mathcal F\subseteq \binom{S}{d}$ be a $d$-uniform set family on a finite ground set $S$.
Define a hypergraph $H(\mathcal F)$ by
$$
V\bigl(H(\mathcal F)\bigr):=\mathcal F,
\qquad
E\bigl(H(\mathcal F)\bigr):=S,
$$
where the edge indexed by $x\in S$ is incident with precisely those vertices $B\in\mathcal F$
for which $x\in B$. Then every vertex of $H(\mathcal F)$ has degree exactly $d$, and
$$
\rk\bigl(H(\mathcal F)\bigr)=\max_{x\in S}\bigl|\{B\in\mathcal F:x\in B\}\bigr|.
$$
Moreover, a $(k+1)$-colouring of the edges of $H(\mathcal F)$ is exactly a $(k+1)$-colouring of
the ground set $S$. In particular, if every $(k+1)$-colouring of $S$ has some $B\in\mathcal F$
and some colour class $C_i$ with
$$
|B\cap C_i|\ge m+1,
$$
then $H(\mathcal F)$ is a $d$-regular hypergraph with no $1/k$-majority $(k+1)$-edge-colouring.
\end{lemma}

\begin{proof}
The first two assertions are immediate from the definition. The final assertion follows because
a vertex corresponding to $B\in\mathcal F$ sees exactly one incident edge for each element of
$B$, so having at least $m+1=\floor{d/k}+1$ incident edges of one colour violates the
$1/k$-majority condition at that vertex.
\end{proof}

\begin{theorem}\label{thm:appendix-lower}
For every fixed integer $k\ge2$, there exists a constant $C_k>0$ such that for all
sufficiently large multiples $d$ of $k$, there exists a $d$-regular hypergraph $H_{k,d}$ with no
$1/k$-majority $(k+1)$-edge-colouring and
$$
\rk(H_{k,d})\le \exp\bigl(dI_k+C_k\log d\bigr).
$$
\end{theorem}

\begin{proof}
Fix a sufficiently large multiple $d=km$ of $k$, and set
$$
t:=d^3,
\qquad
N:=qt=(k+1)t.
$$
For a fixed $t$-subset $T\subseteq [N]$, let
$$
L:=\sum_{j=m+1}^{d}\binom{t}{j}\binom{kt}{d-j}.
$$
Thus $L$ is the number of $d$-subsets of $[N]$ meeting $T$ in at least $m+1$ points.

Choose each member of $\binom{[N]}{d}$ independently with probability
$$
p:=\frac{4N\log q}{L}.
$$
Since $L\ge \binom{t}{m+1}$ and $t=d^3$ while $m+1\le d$, we have
$$
L>4N\log q
$$
for all sufficiently large $d$; in particular, $p<1$. Let $\mathcal F$ denote the resulting random family.

Consider any colouring $c:[N]\to[q]$. By the pigeonhole principle, some colour class of $c$
has size at least $t$; fix a monochromatic $t$-subset $T$. If $\mathcal F$ contains a selected
$d$-set meeting $T$ in at least $m+1$ points, then $c$ cannot correspond to a
$1/k$-majority $(k+1)$-edge-colouring of $H(\mathcal F)$ by Lemma~\ref{lem:incidence-duality}.
For the fixed set $T$, exactly $L$ members of $\binom{[N]}{d}$ have this property, so
$$
\Pr[c\text{ survives}]\le (1-p)^L\le e^{-pL}=q^{-4N}.
$$
Since there are $q^N$ colourings of $[N]$, the union bound gives
$$
\Pr[\exists\text{ surviving }q\text{-colouring of }[N]]\le q^{-3N}<\frac12.
$$

Now fix a point $x\in[N]$. Its degree in $\mathcal F$ has distribution
$$
X_x\sim \Bin\!\left(\binom{N-1}{d-1},p\right)
$$
with mean
$$
\mu:=p\binom{N-1}{d-1}=4d\log q\cdot \frac{\binom{N}{d}}{L}.
$$
Let
$$
\rho:=\frac{L}{\binom{N}{d}}.
$$
We need both upper and lower estimates on $\rho$.

For the upper estimate, let $J$ be the size of the intersection of a uniformly random $d$-set
of $[N]$ with the fixed $t$-set $T$. Then $J$ has hypergeometric distribution with parameters
$(N,t,d)$, mean $d/(k+1)$, and
$$
\Pr[J\ge m+1]=\rho.
$$
A standard comparison theorem of Hoeffding for sampling without replacement
(see~\cite{hoeffding}) gives
$$
\rho\le \Pr\!\left[\Bin\!\left(d,\frac1{k+1}\right)\ge m+1\right]
\le e^{-dI_k}.
$$
Hence
$$
\mu\ge 4d\log q\,e^{dI_k}.
$$
In particular, $\mu\to\infty$ exponentially fast with $d$, so for all sufficiently large $d$,
$$
Ne^{-\mu/3}<\frac12.
$$
By Chernoff's inequality and another union bound,
$$
\Pr[\exists x\in[N]:X_x>2\mu]\le N e^{-\mu/3}<\frac12.
$$

For the lower estimate, note that $L$ contains the term with $j=m+1$, so
$$
\rho\ge \rho_{k,d}:=
\frac{\binom{t}{m+1}\binom{kt}{d-m-1}}{\binom{(k+1)t}{d}}.
$$
Using
$$
\binom{t}{m+1}\binom{kt}{d-m-1}
=
\binom{d}{m+1}\frac{(t)_{m+1}(kt)_{d-m-1}}{d!},
$$
we can rewrite $\rho_{k,d}$ as
$$
\rho_{k,d}=
\binom{d}{m+1}\frac{(t)_{m+1}(kt)_{d-m-1}}{((k+1)t)_d},
$$
where $(x)_s:=x(x-1)\cdots(x-s+1)$. Since $t=d^3$, for $a\in\{1,k,k+1\}$ and every
$0\le s\le d$ we have
$$
\log (at)_s=s\log(at)+O\!\left(\frac{s^2}{t}\right)=s\log(at)+O\!\left(\frac1d\right).
$$
Therefore
\begin{align*}
\log\frac{(t)_{m+1}(kt)_{d-m-1}}{((k+1)t)_d}
&=(m+1)\log\frac1{k+1}+(d-m-1)\log\frac{k}{k+1}+O\!\left(\frac1d\right).
\end{align*}
Let
$$
a_d:=\frac{m+1}{d}=\frac1k+\frac1d.
$$
By Stirling's formula,
$$
-\log \rho_{k,d}
=
d\,D\!\left(a_d\middle\|\frac1{k+1}\right)+O_k(\log d)
=
dI_k+O_k(\log d),
$$
because $a_d=1/k+O(1/d)$ and $D(\cdot\|1/(k+1))$ is smooth at $1/k$. Hence
$$
\rho_{k,d}\ge \exp\bigl(-dI_k-O_k(\log d)\bigr),
$$
and consequently
$$
\mu\le 4d\log q\,\exp\bigl(dI_k+O_k(\log d)\bigr)
=\exp\bigl(dI_k+O_k(\log d)\bigr).
$$

With positive probability, both of the following events occur simultaneously:
\begin{itemize}[leftmargin=1.6em]
\item no colouring $[N]\to[q]$ survives, and
\item every point of $[N]$ belongs to at most $2\mu$ members of $\mathcal F$.
\end{itemize}
Fix such a family $\mathcal F$. By Lemma~\ref{lem:incidence-duality}, the hypergraph
$H(\mathcal F)$ is $d$-regular, has no $1/k$-majority $(k+1)$-edge-colouring, and satisfies
$$
\rk\bigl(H(\mathcal F)\bigr)\le 2\mu
\le \exp\bigl(dI_k+C_k\log d\bigr)
$$
for a suitable constant $C_k>0$.
\end{proof}

\end{document}